\documentclass[11pt]{article}

\usepackage{latexsym}
\usepackage{amsfonts}
\usepackage{amssymb, mathrsfs, amsfonts, amsmath}
\usepackage{amsbsy}
\newcommand\raisepunct[1]{\,\mathpunct{\raisebox{0.5ex}{#1}}}
\newcommand{\diver}{{{\rm{div}\,}}}

\newcommand{\grad}{{{\rm{grad}\,}}}

\newcommand{\Hess}{{{\rm Hess}\,}}

\newtheorem{theorem}{\bf Theorem}[section]
\newtheorem{lemma}[theorem]{\bf Lemma}
\newtheorem{corollary}[theorem]{\bf Corollary}
\newtheorem{proposition}[theorem]{\bf Proposition}

\begin{document}

\title{The spectrum of the Laplacian and volume growth of proper minimal submanifolds}

\author{G. Pacelli Bessa\thanks{\small Partially supported by CNPq-Brazil grant\# 	303057/2018-1.}  \and Vicent Gimeno\thanks{\small Work  partially supported by the Research Program of University Jaume
I Project P1-1B2012-18, and DGI -MINECO grant (FEDER) MTM2013-48371-C2-2-P.} \and Panagiotis Polymerakis\thanks{\small Supported by the Max Planck Institute for Mathematics in Bonn.}}
\date{\today}
\maketitle
\begin{abstract}We give  upper bounds for the bottom of the essential spectrum of properly immersed minimal submanifolds of $\mathbb{R}^{n}$ in terms of their volume growth. Our result  improves the extrinsic version of  Brook's  essential spectrum estimate given by Ilias-Nelli-Soret in \cite[Cor.3]{ins}.
\vspace{.1cm}


\vspace{.1cm}

\end{abstract}
\section{Introduction}
Let $M$ be a complete Riemannian $n$-manifold and let $\triangle = \diver \circ \grad$ be the Laplace-Beltrami operator (Laplacian) acting on $C_{0}^{\infty}(M)$ the space of smooth functions with compact support. The Laplacian has a unique self-adjoint extension to an operator $\triangle \colon \mathcal{D}(\triangle)\to L^{2}(M)$ whose domain is  $\mathcal{D}(\triangle)=\{f\in L^{2}(M)\colon  \triangle f \in L^{2}(M)\}$. The spectrum  $\triangle$ is the set of  $\lambda\in \mathbb{R}$ for which  ${\rm Ker}(\triangle +\lambda I)\neq \{0\}$  or  $(\triangle +\lambda I)^{-1}$ is unbounded. We will  refer to $\sigma(\triangle)$ as the spectrum of $M$ and denote it by $\sigma(M)$.   Those $\lambda$'s for which ${\rm Ker}(\triangle +\lambda I)\neq \{0\}$ are the eigenvalues of $M$ and the elements of ${\rm Ker}(\triangle +\lambda I)$ are the eigenfunctions associated to $\lambda$.
The set of all eigenvalues of $M$  is the point spectrum $\sigma_{p}(M)$ and the subset of the point spectrum formed by the isolated eigenvalues  with finite multiplicity (${\rm dim \,Ker}(\triangle +\lambda I)<\infty$) is called the discrete spectrum $\sigma_{d}(M)$. The essential spectrum is  $\sigma_{ess}(N)=\sigma(M)\setminus \sigma_{d}(M)$, see \cite{davies}. 
 When $M$ is compact the spectrum of $\triangle$ is discrete while when $M$ is non-compact the spectrum may be purely continuous,  ($\sigma_{p}(M)=\emptyset$), like the Euclidean space  $\mathbb{R}^{n}_{\raisepunct{,}}$ purely  discrete, ($\sigma_{ess}(M)=\emptyset$), as the simply connected Riemannian manifolds with highly negative curvature, \cite{donnelly-li}  or  may be a mixture of both types, \cite{donnelly81a, donnelly90}. The very basic question \cite{schoen-yau} is {\em for what geometries $\inf \sigma(M)>0$? } It was shown by
McKean in \cite{mckean}  that if $M$ is a simply connected Riemannian manifold with curvature $K_{M}\leq -\delta^{2}<0$ then $\inf \sigma (M)\geq (n-1)^2\delta^2/4$.  Cheng   has  shown in \cite{cheng75} that if $M$ is complete with non-negative Ricci curvature ${\rm Ric}_{M}\geq 0$ then $\inf \sigma (M)=0$.  An curvature free estimate for the bottom of the spectrum was obtained by R. Brooks in \cite{brooks}. More precisely, 
let $M$ be a complete Riemannian manifold of infinite volume and
   $v(r)={\rm vol}(B_p(r))$ be the volume of  the geodesic ball  $B_p(r)$ of radius $r$ centred at $p\in M$. Set $$\mu =\limsup_{r\to \infty}\frac{\log(v(r))}{r}\cdot$$    Brooks   proved  that $\inf \sigma_{ess}(M)\leq \mu^{2}/4$. 
   
   It is a classical result  due to   Efimov-Hilbert    that complete surfaces with curvature $K_{_M}\leq -\delta^2<0$ can not be isometrically immersed in $\mathbb{R}^{3}$, see \cite{efimov, hilbert}. Naturally one is lead to ask if  a complete surface with $\inf \sigma (M)>0$ can be isometrically immersed in $\mathbb{R}^{3}_{\raisepunct{.}}$ 
   It  turns out that the  examples,   constructed by Nadirashvili \cite{nadirashvili} and the Spanish School of Geometry  \cite{alarcon,lopez-martin-morales,lopez-martin-morales2,pacomartin1,pacomartin2}, of bounded complete minimal surfaces of $\mathbb{R}^{3}$
 have $\inf \sigma(M)>0$, see   \cite{bessa-montenegro, bessa-montenegro2007,silvana}. However, the question whether a manifold with $\inf \sigma(M)>0$ can be minimally and properly immersed in the Euclidean space can is still valid. In some sense, it complements the question  raised by  S. T. Yau in \cite[p.240]{yau}  when he asked what upper bounds can one give for the bottom of the spectrum of complete  immersed minimal surfaces in $\mathbb{R}^{3}_{\raisepunct{.}}$  S. Ilias, B. Nelli and M. Soret \cite[Cor.3]{ins} gave a fairly general answer to Yau's question establishing  a Brook's type upper estimate for the bottom of the essential  spectrum of any properly immersed submanifold  $\varphi\colon M \to \mathbb{R}^{n}_{\raisepunct{,}}$  with infinite volume. They proved that\[ \inf \sigma_{ess}(M)\leq \left[\liminf_{r\to \infty} r^{-1}\log({\rm vol} (\Omega_r))\right]^2/4.\] Here $\Omega_r=\varphi^{-1}(B_p^{^{N}}(r))$ of radius $r$. 
In this article we prove a    strong Brook's type upper estimate for the bottom of the essential spectrum of  properly immersed minimal $m$-submanifolds  of the Euclidean $n$-space. Indeed, letting $\Omega_r\subset M$ be the  extrinsic geodesic ball  of radius $r>0$  of  a  properly minimal immersion $\xi \colon M \hookrightarrow \mathbb{R}^{n}$ of a complete Riemannian $m$-manifold $M$ into  $\mathbb{R}^{n}$    with $\xi (p)=o$, i.e. $\Omega_r=\xi^{-1}(B_o(r))$,  we prove the following  result.
\begin{theorem}\label{thm1}
Let $\xi\! \colon\! M \!\to \mathbb{R}^{n}$ be a proper  isometric  minimal immersion of a complete  $m$-submanifold $M$ of  $\mathbb{R}^{n}$ with $\xi (p)=o$. The bottom of the essential spectrum  is bounded above by
\[\inf \sigma_{ess}(M)\leq m\cdot\liminf_{r\to \infty}\left[r^{-2}  \log ({\rm vol}(\Omega_r))\right].\]
\end{theorem}
Theorem \ref{thm1} has a number of corollaries.
Let $\Theta(r)={\rm vol}(\Omega_r)/{\rm vol}(B^{m}(r),$ where $B^{m}(r)$ is a geodesic ball of radius $r$ in the Euclidean space $\mathbb{R}^m_{\raisepunct{.}}$ In \cite{lima-mari-montenegro-viera}, Lima et al.,  proved that if $\liminf_{s\to \infty}(\log(\Theta(s))/\log(s))=0$ then  $\sigma(M)=[0, \infty)$. In particular $\inf\sigma_{ess}(M)=0$.  
\begin {corollary}Let $\xi\! \colon\! M \!\to\! \mathbb{R}^{n}$ be a complete properly immersed  minimal     $m$-submanifold $M$ of $\mathbb{R}^{n}$ with $\xi (p)=o$. If \[\liminf_{r\to \infty}\frac{\log(\Theta(r))}{r^2}=0\] then $\inf\sigma_{ess}(M)=0$.\label{cor1}
\end{corollary}

\begin{corollary}Let $\xi\! \colon\! M \!\to\! \mathbb{R}^{3}$ be a complete properly embedded  minimal     surface $M$ of $\mathbb{R}^{3}$ with $\xi (p)=o$ and  let $\kappa(r)=\inf_{x\in\Omega_r}\left\{  K_{_M}(x)\right\}$, where  $K_{_M}(x)$  is the Gaussian curvature of $M$ at $x$. If  \[\liminf_{r\to \infty}\frac{\log(\vert\kappa(r)\vert)}{r^2}=  0\]  then $\inf\sigma_{ess}(M)=0$.\label{cor2}
\end{corollary}

\begin{corollary} Let  $\xi\colon M \to \mathbb{R}^n$ be an isometric minimal immersion of   a complete  $m$-dimensional Riemannian manifold into $\mathbb{R}^{n}$ with $\xi(p)=o$. Suppose  that for some $\sigma>0$ \[\int_{M}e^{-\sigma \Vert \xi (x)\Vert^2}d\mu (x)<\infty.\] Then the immersion $\xi$ is proper, see \cite[Thm.1.1]{GP}, and  $$\inf\sigma_{ess}(M)\leq m\sigma.$$\label{cor3}
\end{corollary}

\section{Poof of Theorem  \ref{thm1}}A  {\em model  $n$-manifold} $\mathbb{M}_{h}^{n}$, with radial sectional curvature $-G(r)$ along the geodesics  issuing from the origin, where $G\colon \mathbb{R}\to \mathbb{R}$ is a smooth even function, is the quotient space
\[
\mathbb{M}^n_{h}=[0,R_h)\times\mathbb{S}^{n-1}/\sim
\] with $(\rho, \theta)\sim (\tilde{\rho},\beta)\Leftrightarrow \rho=\tilde{\rho}=0$  or $ \rho=\tilde{\rho}$ and $\theta=\beta$, endowed with the metric $ds^2_{h}=d\rho^2+h^2(\rho)d\theta^2$
%where $\sigma$ denotes the unique solution of the Cauchy problem
%\begin{eqnarray}\label{ubs2}
%\left\{\begin{array}{l}
%h''-Gh=0,\\
%h(0)=0, h'(0)=1
%\end{array}\right.
%\end{eqnarray}on $[0,R_h)$
where  $h\colon [0, \infty) \to \mathbb{R}$ is the unique solution of the Cauchy problem \begin{eqnarray*}%\label{ubs2}
\left\{\begin{array}{rll}
h''-Gh&=&0,\\
h'(0)&=&1,\\ h^{(2k)}(0)&=&0,\,\, k=0,1,\ldots,
\end{array}\right.
\end{eqnarray*} and
$R_h$ is the largest positive real number such that $h\vert_{(0,R_{h})}>0$. The model $\mathbb{M}_{h}^{n}$ is non-compact with pole at the origin  $o=\{0\}\times \mathbb{S}^{n-1}/\sim$ if
 $R_h=\infty$. Observe that $\mathbb{M}^{n}_{t}=\mathbb{R}^{n}$,
  $\mathbb{M}^{n}_{\sinh(t)}=\mathbb{H}^{n}(-1)$ and if $h(t)=\sin(t)$ and $R_h=\pi$ then $\mathbb{M}^{n}_{\sin(t)}=\mathbb{S}^{n}$.
 
 \vspace{2mm}
 
  If $G$ satisfies \[\begin{array}{lll}G_{-}\in L^{1}(\mathbb{R}^{+})&{\rm and}& \int_{t}^{\infty}G_{-}(s)ds \leq \displaystyle\frac{1}{4t}\end{array}\] then $h'>0$ in $\mathbb{R}^{+}$ and $\mathbb{M}_{h}^{n}$ is geodesically complete, \cite[Proposition 1.21]{BMR}.
 The geodesic ball centered at the origin  with radius $r<R_h$ is the set $B_h(r)=[0,r)\times \mathbb{S}^{n-1}/\sim$ whose volume  and the volume of its boundary are given respectively, by \[ \begin{array}{lll}
V(r)=\omega_{n}\int_{0}^{r}h^{n-1}(s)ds &  {\rm and} &
S(r)= \omega_{n}h^{n-1}(r),
\end{array} \] where $\omega_n={\rm vol}(\mathbb{S}^{n-1})$.
The  Laplace operator on $B_h(r)$, expressed in polar coordinates $(\rho, \theta)$, is given by \[\triangle= \frac{\partial^{2}}{\partial \rho^{2}}+ (n-1)\frac{h'}{h}\frac{\partial }{\partial \rho} + \frac{1}{h^{2}}\triangle_{\theta}.\]
Let $\rho(x)={\rm dist}_{\mathbb{M}_{h}^{n}}(o,x)$ be the distance function to  the origin $o$ on $\mathbb{M}^{n}_{h}$. The hessian of $\rho$ is given by the following expression
\begin{equation}\Hess \rho (x)(e_i,e_i)=\frac{h'}{h}(\rho(x))\left\{\langle e_i,e_i\rangle - d\rho\otimes d\rho (e_i,e_i)\right\}, \label{eq3a}\end{equation} where $\{ e_1, \ldots, e_m\}$ is an orthormal basis of $T_x\mathbb{M}_{h}^{n}$. Let 
$\varphi \colon M \hookrightarrow \mathbb{M}_h^{n}$ be an isometric immersion of a complete    $m$-manifold into $\mathbb{M}_h^{n}$. Suppose that $\varphi(p)=o$ for some  $p\in M$. The function $t\colon M \to \mathbb{R}$ given by $t(y)=\rho \circ \varphi (y)$ is smooth in $M\setminus \varphi^{-1}(o)$. The hessian of $t$ is given by
\[ \Hess_{_M}t(q)(e_i,e_i)=\Hess_{_{\mathbb{M}^n_{h}}}\rho( e_i, e_i)+\langle \grad \rho, \alpha (e_i,e_i)\rangle.\] Here we are identifying $e_i= d\varphi \cdot e_i$, see \cite{JK}. 
In particular, \begin{equation}\triangle_{_M}t (q)=\sum_{i=1}^{m}\Hess_{_{\mathbb{M}^n_{h}}}\rho(d\varphi \cdot e_i, d\varphi \cdot e_i)+\langle \grad \rho, \stackrel{\rightarrow}{H} \rangle. \label{eq4}\end{equation}
Let $ \varphi\colon M \to \mathbb{M}^{n}_{h}$ be a complete properly and minimally immersed $m$-submanifold of  $\mathbb{M}_{h}^{n}$ with radial sectional  curvature $-G(\rho)\leq 0$, $\varphi(p)=o$ and let $\Omega_r$ be the pre-image $\varphi^{-1}(B_o(r))$.

\begin{lemma}\label{lemma1}For almost any $r>0$ we have that
\begin{equation*}\int_{\partial \Omega_r}\vert \grad t\vert\, d\nu\leq 
m\, \frac{ h'}{h}(r)\,{\rm vol}(\Omega_r).\end{equation*}\end{lemma}
Proof:
 Let $\phi\colon M \to \mathbb{R}$ given by \[\phi(y)=\left(\displaystyle\int_{0}^{t}h(s)ds\right)\circ \rho \circ \varphi(y)=\displaystyle\int_{0}^{\rho(\varphi(y))}h(s)ds.\]
At a point $q\in M$ and an orthonormal basis $\{e_1,\ldots, e_m\}$ of $T_qM$ that, using \eqref{eq3a} and \eqref{eq4}, we have 
 \begin{eqnarray}\triangle_{_M}\phi &=&\sum_{i=1}^{m}\left[\phi''(\rho)\langle \grad \rho, e_i\rangle^2 +\phi'(\rho)\frac{h'}{h}(\rho)\left\{\langle e_i, e_i\rangle - \langle \grad \rho, e_i\rangle^2\right\}\right]\nonumber \\
 &&\nonumber \\
 &=&m\, h'(\rho).\nonumber 
\end{eqnarray}Since $G\geq 0$, we have that $h''(s)=G(s)h(s) \geq 0$ for $s>0$, which implies that $h'$ is non-decreasing. In view of Sard's theorem, $\Omega_r$ is smoothly bounded for almost any $r>0$. For any such $r$, we compute 
 \begin{eqnarray*}m\,h'(r) {\rm vol}(\Omega_r) &\geq & \int_{\Omega_r}\triangle_{_M}\phi\, d\mu  \nonumber \\
 &=& \int_{\partial \Omega_r}\left\langle \grad \phi, \frac{\grad t}{\vert \grad t \vert}\right\rangle d\nu\\
 &=&h(r) \int_{\partial \Omega_r}\vert   \grad t \vert d\nu.\nonumber 
\end{eqnarray*}Thus \[  \int_{\partial \Omega_r}\vert   \grad t \vert d\nu\leq m\frac{h'(r)}{h(r)}{\rm vol}(\Omega_r). \]This proves Lemma \ref{lemma1}.

\vspace{2mm}

Let $\mathcal{H}_0^1(M)$ be the space of square-integrable functions with square-integrable gradient. Given a non-zero $u \in \mathcal{H}_0^1(M)$, set
\[
\mathcal{R}(u) = \frac{\int_M \| \grad  u \|^2}{\int_M u^2}.
\]
For $r>0$, define $u_r \colon M \to [0,+\infty)$ by
\[
u_r(x) = \begin{cases}
	\phi (r)-\phi(t(x)) & \text{ for } x \in \Omega_r, \\
	0 & \text{else}.
\end{cases}
\]
It should be noticed that $u_r \in \mathcal{H}_0^1(M)$, being compactly supported and Lipschitz. Consider also the function
\[
v_r = \frac{u_r}{\left(\int_M u_r^2\right)^{1/2}}.
\]
This renormalization gives rise to sequences of functions which converge weakly to zero as the next lemma indicates.
\begin{lemma}\label{weak convergence}
For any sequence $(r_n)_{n \in \mathbb{N}} \subset (0,+\infty)$ with $r_n \rightarrow + \infty$ we have that $v_{r_n} \rightharpoonup 0$ in $L^2(M)$.
\end{lemma}
Proof: For any $c > 0$ we compute
\[
\int_{\Omega_c} v_{r_n}^2 = \frac{\int_{\Omega_c} (r_n^2 - t^2)^2}{\int_{\Omega_{r_n}} (r_n^2 - t^2)^2} \leq \frac{r_n^4 {\rm vol}(\Omega_c)}{\int_{\Omega_{r_n/2}} (r_n^2 - t^2)^2} \leq \frac{16  {\rm vol}(\Omega_c)}{9  {\rm vol}(\Omega_{r_n/2})} \rightarrow 0.
\]
Keeping in mind that $(v_{r_n})_{n \in \mathbb{N}}$ is bounded in $L^2(M)$, this shows that $v_{r_n} \rightharpoonup 0$. This completes the proof Lemma \ref{weak convergence}.

\vspace{2mm}

The significance of considering sequences that converge weakly to zero in order to estimate the bottom of the essential spectrum is illustrated in the following.

\begin{proposition}\label{estimate}
Consider $(v_n)_{n \in \mathbb{N}} \subset \mathcal{H}_0^1(M)$ with $\| v_n \|_{L^2(M)} = 1$ and $v_n \rightharpoonup 0$ in $L^2(M)$. Then the minimum of the essential spectrum of $M$ is bounded by
\[
\inf\sigma_{ess}(M) \leq \liminf_{n} \mathcal{R}(v_n).
\]
\end{proposition}
Proof: If the right hand side is infinite, there is nothing to prove. If it is finite, we denote it by $\lambda$ and after passing to a subsequence if necessary, we may suppose that $\mathcal{R}(v_n) \rightarrow \lambda$. Assume to the contrary that $\inf \sigma_{ess}(M) > \lambda$. Then
\[
\sigma(\Delta) \cap (- \infty , \lambda] = \sigma_d(\Delta) \cap [0 , \lambda] = \{\lambda_1 , \dots, \lambda_k\},
\]
where $\lambda_i$'s are eigenvalues of the unique self-adjoint extension $-\triangle$ of minus the Laplacian of finite multiplicity, for some $k \in \mathbb{N}$. Let $E_i$ be the eigenspace corresponding to $\lambda_i$ and denote by $E$ their sum. Then the spectrum of the restriction $\triangle|_{E^\perp}$ of $\triangle$ to the $L^2$-orthogonal complement of $E$ is given by $\sigma(\triangle|_{E^\perp}) = \sigma(\triangle) \smallsetminus \{\lambda_1 , \dots, \lambda_k\}$ and in particular, the minimum of its spectrum is greater than $\lambda$.
	
Writing $v_n = u_n + w_n$ with $u_n \in E$ and $w_n \perp E$, we readily see that $u_n \rightarrow 0$ and $\triangle u_n \rightarrow 0$ in $L^2(M)$, since $v_n \rightharpoonup 0$ and $E$ is finite dimensional. This implies that 
	\[
	\int_M |\grad u_n|^2 = -\int_M u_n \triangle u_n \rightarrow 0.
	\]
	Moreover, we obtain that $\| w_n \|_{L^2(M)} \rightarrow 1$ and
	\[
	\int_{M} | \grad w_n|^2 = -\int_M | \grad (v_n - u_n)|^2 \rightarrow \lambda.
	\]
	We conclude that $\mathcal{R}(w_n) \rightarrow \lambda$, which yields that the minimum of the spectrum of $\triangle|_{E^\perp}$ is less or equal to $\lambda$, which is a contradiction, that establishes Proposition \ref{estimate}.
	
	\vspace{2mm}

Our goal now is to estimate $\mathcal{R}(v_r) = \mathcal{R}(u_r)$. To this end, using the co-area formula and Lemma \ref{lemma1}, we compute
%Define $u\colon \Omega_r\to [0,\infty)$ given by $u(x)=\phi (r)-\phi(t(x))$. We have \[\lambda_1(\Omega_r)\leq \frac{\int_{\Omega_r}\vert \grad u\vert^2\,d\mu}{\int_{\Omega_r}u^2d\mu}\]On the other hand, by the co-area formula
\begin{eqnarray}
\int_{\Omega_r}\vert \grad u\vert^2\,d\mu&=& \int_{\Omega_r}h^2(t(x))\vert \grad t \vert^2 d\mu\nonumber \\
&=& \int_{0}^{r}h^2(s)\int_{\partial \Omega_s}\vert \grad t \vert \,d\nu ds\nonumber \\
&\leq & m \int_{0}^{r}h(s)h'(s){\rm vol}(\Omega_s)d s\nonumber\\
&\leq & m\, h'(r)\int_{0}^{r}h(s){\rm vol}(\Omega_s)d s.\nonumber 
\end{eqnarray} It follows from \cite[Lemma 2, Propositions 2 and 3]{lima-mari-montenegro-viera} that ${\rm vol}(\Omega_s)$ is locally absolutely continuous with  \[\frac{d\, {\rm vol}(\Omega_s) }{ds}=\int_{\partial \Omega_s}\frac{1}{\vert \grad t\vert }d\nu,\]
and
\begin{eqnarray}
\int_{\Omega_r}u^2d\mu&=&\int_{0}^{r}\left(\phi(r)-\phi(s)\right)^2\int_{\partial \Omega_s}\frac{1}{\vert \grad t\vert }d\nu ds\nonumber \\
&=& \int_{0}^{r}\left(\phi(r)-\phi(s)\right)^2\,\frac{d \,{\rm vol}(\Omega_s)}{ds}ds\nonumber\\
&=& 2\int_{0}^{r}(\phi(r)-\phi(s))\phi'(s)\,{\rm vol}(\Omega_s)ds.\nonumber 
\end{eqnarray}
Thus,
\[\mathcal{R}(u_r)\leq \frac{ m\,h'(r) \displaystyle\int_{0}^{r}h(s)\,{\rm vol}(\Omega_s)d s}{2\displaystyle\int_{0}^{r}(\phi(r)-\phi(s))h(s)\,{\rm vol}(\Omega_s)ds}\cdot\]
Letting $F(r)=\displaystyle\int_{0}^{r}(\phi(r)-\phi(s))h(s)\,{\rm vol}(\Omega_s)ds$ we have that \[F'(r)= h(r)\int_{0}^{r}h(s)\,{\rm vol}(\Omega_s)ds.\] Therefore,
\[\mathcal{R}(v_r) = \mathcal{R}(u_r)\leq \frac{m}{2}\frac{h'(r)}{h(r)}\frac{F'(r)}{F(r)}\cdot\]
When  $h(s)=s$, i.e. the model $\mathbb{M}^{m}_{h}=\mathbb{R}^{m}_{\raisepunct{,}}$ this inequality reads as
\[
\mathcal{R}(v_r) \leq m \frac{(\log F(r))^\prime}{2r}
\]
for any $r > 0$. We deduce from Lemma \ref{weak convergence} and Proposition \ref{estimate} that
\begin{equation}\label{estimate with derivatives}
\inf \sigma_{ess}(M) \leq m \liminf_{r \rightarrow + \infty} \frac{(\log F(r))^\prime}{2r}\cdot
\end{equation}

Consider any $c \in \mathbb{R}$ with
\[
c \leq \liminf_{r \rightarrow + \infty}  \frac{(\log F(r))^\prime}{2r}\cdot
\]
Then for any $\varepsilon > 0$ there exists $r_0 > 0$ such that
\[
(\log F(r))^\prime \geq 2(c - \varepsilon) r
\]
for any $r \geq r_0$. Integrating gives that
\[
\log F(r) - \log F(r_0) \geq (c - \varepsilon) (r^2 - r_0^2)
\]
for any $r \geq r_0$, which yields that
\[
\liminf_{r \rightarrow + \infty}  \frac{\log F(r)}{r^2} = \liminf_{r \rightarrow + \infty}  \frac{\log F(r) - \log F(r_0)}{r^2 - r_0^2} \geq c - \varepsilon.
\]
We conclude from this together with (\ref{estimate with derivatives}) that
\[
\inf \sigma_{ess}(M) \leq m \liminf_{r \rightarrow + \infty} \frac{\log F(r)}{r^2}\cdot
\]
This establishes Theorem \ref{thm1} after noticing that $F(r) \leq {\rm vol}(\Omega_r)\,r^4/8$.
%This implies that\[\frac{2\lambda_1(\Omega_r)}{m}\frac{h(r)}{h'(r)} \leq \log (F(r))'.\] Choose $r_o$ so that $F(r_o)=1$ and integrate from $r_o$ to $r$ we obtain that
% \[\displaystyle\frac{2\lambda_1(\Omega_r)}{m}\int_{r_o}^{r}\frac{h(s)}{h'(s)}ds\leq \log(F(r))\]On the other hand \begin{eqnarray}F(r)&\leq & {\rm vol}(\Omega_r)\int_{0}^{r}(\phi(r)-\phi(s))h(s)\, ds\nonumber  \end{eqnarray}
%Thus \begin{equation} \displaystyle\frac{2\lambda_1(\Omega_r)}{m}\int_{r_o}^{r}\frac{h(s)}{h'(s)}ds\leq \log\left( {\rm vol}(\Omega_r)\int_{0}^{r}(\phi(r)-\phi(s))h(s)\, ds\right)\label{eq7}
%\end{equation}This yields 
%\[\lambda^{\ast}(\Omega_r)\leq \frac{m \log\left( {\rm vol}(\Omega_r)\int_{0}^{r}(\phi(r)-\phi(s))h(s)\, ds\right)}{2\int_{r_o}^{r}\frac{h(s)}{h'(s)}ds}\cdot\]
%When  $h(s)=s$, i.e. the model $\mathbb{M}^{m}_{h}=\mathbb{R}^{m}_{\raisepunct{,}}$, the inequality \eqref{eq7} reads as
%\begin{equation}\label{eq8} \lambda^{\ast}(\Omega)\leq  \frac{m\log\left(r^4 {\rm vol}(\Omega_r)/8\right)}{r^2-r_o^2}\end{equation} Observe that \[ 
%\liminf_{r\to \infty} \left(\frac{\log( {\rm vol}(\Omega_r))}{r^2}\right)= \liminf_{r\to \infty} \frac{\log\left(r^4 {\rm vol}(\Omega_r)/8\right)}{r^2-r_o^2}.\] Since $\lim_{r\to \infty}\lambda^{\ast}(\Omega_r)=\lambda^{\ast}(M)$  we have that \[\lambda^{\ast}(M)\leq m\liminf_{r\to \infty} \left(\frac{\log( {\rm vol}(\Omega_r))}{r^2}\right)\]
%This proves the  theorem. 
\vspace{2mm}

To prove the corollaries we proceed as follows. Observe that \[{\rm vol}(\Omega_r)=\Theta(r){\rm vol}(B^{m}(r))\] then $$\frac{\log({\rm vol}(\Omega_r))}{r^2}=\frac{\log(\Theta(r))}{r^2}+\frac{\log({\rm vol}(B^{m}(r))}{r^2}\cdot$$Thus by Theorem \ref{thm1} $$\inf\sigma_{ess}(M)\leq m\liminf_{r\to \infty}\left(\frac{\log(\Theta(r))}{r^2}\right).$$ This proves Corollary \ref{cor1}. 
\vspace{2mm}

Given a unit normal vector field $N:M\to\mathbb{R}^3$, consider the  tubular neighbourhood of $\varphi(\Omega_r)$,
 $$T_{\epsilon}(\Omega_r)=\{y\in \mathbb{R}^{3}\, \colon\, y=q+x N(q),\quad -\epsilon<x<\epsilon,\quad q\in M\}.$$    By \cite[p.9]{Gray} the volume of $T_{\epsilon}(\Omega_r)$ for $\epsilon$ small enough is given by  \begin{eqnarray}{\rm vol}(T_\epsilon(\Omega_r))&=&2\epsilon\,{\rm vol}(\Omega_r)+ \frac{2\epsilon^2}{3}\int_{\Omega_r}K_{_M}(x)d\mu(x)\nonumber \\
 &\leq &{\rm vol}(B^{3}(\epsilon+r))\nonumber \\
 &=& \frac{4\pi}{3}\left(\epsilon+r\right)^3.\nonumber\end{eqnarray}On the other hand 
 \begin{eqnarray}
 2\epsilon\,{\rm vol}(\Omega_r)+ \frac{2\epsilon^2}{3}\int_{\Omega_r}K_{_M}(x)d\mu(x)\geq 2\epsilon\left(1+\frac{\epsilon}{3}\kappa(r)\right){\rm vol}(\Omega_r).\nonumber  
 \end{eqnarray}Choosing $\epsilon \leq -\frac{3\alpha}{\kappa(r)}$, $0<\alpha<1$ we have that 
 \begin{eqnarray}
 {\rm vol}(\Omega_r)&\leq & \frac{4\pi}{3}\frac{\left(r-\frac{3\alpha}{\kappa(r)}\right)^3}{-\frac{6\alpha}{\kappa(r)}(1-\alpha)}\nonumber\\
 &=& \frac{4\pi}{18\alpha (1-\alpha)}\left(r-\frac{3\alpha}{\kappa(r)}\right)^3(-\kappa(r)).\nonumber
 \end{eqnarray}Thus
 \begin{equation*}
 \liminf_{r\to \infty}\frac{\log({\rm vol}(\Omega_r))}{r^{2}}\leq \liminf_{r\to \infty}\frac{\log(\vert \kappa(r)\vert)}{r^{2}}\cdot
 \end{equation*}This proves Corollary \ref{cor2}.
 \vspace{3mm}
 
Suppose that $C = \int_{M}e^{-\sigma \Vert \xi (x)\Vert^2}d\mu (x)<\infty$ for some $\sigma>0$. Then
\[
C \geq \int_{\Omega_r} e^{-\sigma \Vert \xi (x)\Vert^2}d\mu (x) \geq e^{-\sigma r^2} {\rm vol}(\Omega_r)
\]
for any $r>0$. We derive from \cite[Thm. 1.1]{GP} that the immersion $\xi$ is proper, and the proof of Corollary \ref{cor3} is completed by Theorem \ref{thm1}.

\vspace{0.3cm}
{\small
\begin{tabular}{l p{0.5\linewidth}}
G. Pacelli Bessa  & Vicent Gimeno\\
Departamento de Matem\'atica  &Departament de Matem\`atiques- IMAC\\
Universidade Federal do Cear\'a & Universitat Jaume I                           \\
60455-760-Fortaleza, Brazil& Castell\'{o}, Spain.                         \\
e-mail: bessa@mat.ufc.br & e-mail:  gimenov@uji.es\\
\\
\\
Panagiotis Polymerakis & \\
Max Planck Insitute for Mathematics& \\
Vivatsgasse 7, 53111, Bonn. & \\
e-mail: polymerp@mpim-bonn.mpg.de & 
\end{tabular}
}

\end{document}